\documentclass[11pt,a4paper]{article}
\usepackage[latin1]{inputenc}
\usepackage{amsmath}
\usepackage{amsfonts}
\usepackage{amssymb}
\usepackage{graphicx}
\usepackage[top=2cm, bottom=2cm, left=2cm, right=2cm]{geometry}
\begin{document}
\thispagestyle{empty}	
\textbf{{\Large Note On The Maximal Prime Gaps}}\\
N. A. Carella
\vskip .5 in

\textbf{\textit{Abstract}:} This note presents a result for the prime gap of the form \(p_{n+1}-p_n\leq c \log \left(p_n\right){}^{1+\epsilon
}\), where \(c>0\) is a constant, for any arbitrarily small real number \(\epsilon >0\), { }and all sufficiently large integer \(n\geq n_0(\epsilon
)\). Equivalently, the result shows that short intervals \([x,x+y]\) contain prime numbers for all sufficiently large real numbers \(x\geq x_0(\epsilon
)\), and \(y\geq c \log (x)^{1+\epsilon }\) unconditionally. An application demonstrates that a prime \(p\geq x\geq 2\) can be determined in deterministic
polynomial time \(O\left(\log (x)^8\right)\) . \\

\textit{Mathematics Subjects Classification}: 11A41, 11Y55, 11N05, 11N13, 11P32.

\vskip 1 in
\section{Introduction}
Let \(n\geq 2\) be an integer, and let \(p_n\geq 2\) denotes the \textit{ n}th prime number. The \textit{n}th prime gap is defined as the difference
\(d_n=p_{n+1}-p_n\) of two consecutive prime numbers. This note investigates the maximal size of the prime gap function. For a large real number
\(x\geq 1\), the maximal prime gap function is defined by\\ 
\begin{equation}
d_{\max }=\max  \left\{ p_{n+1}-p_n\leq x \right\}.
\end{equation} \\
The average prime gap for a short sequence of consecutive primes \(2,3,5,7, \text{...},p_n,p_{n+1}\leq x\) is given by the asymptotic formula\\ 
\begin{equation}
\frac{x}{\pi
(x)}=\frac{1}{x}\sum _{p_n\leq x} \left(p_{n+1}-p_n\right) =\log  x +O(\log \log x).
\end{equation} \\
This immediately shows that \(d_{\max }\geq \log  x +O(\text{loglog} x)\). Some of the leading prime gaps open problems provide upper bounds and
lower bounds of the maximal prime gap.\\

\textbf{Cramer Prime Gap Problem.} This problem claims that \\
\begin{equation}
p_{n+1}-p_n\leq c_0\log \left(p_n\right)^2
\end{equation}\\
for some constant \(c_0>0\), see \cite{CH37}. The precise implied constant in this formula has been debated in the literature, this conjecture is discussed
in \cite{GR95}.\\

\textbf{Erdos-Rankin Prime Gap Problem.} This problem claims that \\
\begin{equation}
p_{n+1}-p_n\geq c_1\log \left(p_n\right)\log _2\left(p_n\right)\log _4\left(p_n\right)\log
_3\left(p_n\right)^{-2}
\end{equation}\\
for some constant \(c_1>0\), and that the constant cannot be replaced with any increasing function of \(n\geq 1\), see \cite{EP35}, \cite{EP81}, \cite{MH85}, \cite{RR38}. The notation \(\log _k(x)=\log (\text{$\cdots $log}(x)\cdots )\) denotes the k-fold iterations of the logarithm.\\

The Cramer conjecture provides an upper bound of the maximal prime gap function \(d_{\max }\leq c_0\log \left(p_n\right){}^2\). In comparison,
the Erdos conjecture provides a lower bound of the maximal size of the prime gap function \(d_{\max }\geq c_1\log \left(p_n\right)\log _2\left(p_n\right)\log
_4\left(p_n\right)\log _3\left(p_n\right){}^{-2}\). Extensive details on the theory of maximal prime gap and primes in short intervals are given
in \cite{RP96}, \cite{BP01}, \cite{HM72}, \cite{SK07}, \cite{SA43}, \cite{GP76}, \cite{GR95}, \cite{GP07}, \cite{OR99}, et alii. Some information on the calculations of { }the ratio \(\lim_{n\longrightarrow \infty
} \left.d_n\right/\log  p_n\), called the merit factor,  appears in [\cite{NW00}, p. 340], and \cite{GP07}.  \\

\textbf{Note:}  A pair of new results, proved by independent authors, have disproved the Erdos-Rankin prime gap problem, see \cite{FK14} and \cite{MJ14}.\\

\section{The Proof of the Result}
The numerical evidence, \cite{NN03}, \cite{PP16}, and the result within suggest that the Erdos conjecture is closer to the correct order of magnitude of the maximal
prime gap function:\\
\begin{equation}
c_1\frac{\log \left(p_n\right)\log _2\left(p_n\right)\log _4\left(p_n\right)}{\log _3\left(p_n\right)^2}\leq d_{\max
}\leq \log \left(p_n\right)^{1+\epsilon }
\end{equation}\\
as seemed to be confirmed below.\\

\textbf{Theorem 1.} \textit{Let \(\epsilon >0\) be an arbitrarily small real number, and let \(n\geq 1\) be an integer. Then, the nth prime gap satisfies}\\ 
\begin{equation}
p_{n+1}-p_n\leq c \log \left(p_n\right)^{1+\epsilon }
\end{equation}\\
\textit{where \(c>0\) is a constant, for all sufficiently large integer \(n\geq n_0(\epsilon )\) unconditionally.}\\

\textbf{Proof:} Given a small real number \(\epsilon >0\), let \(n\in \mathbb{N}\) be an integer, and let \(p_{n+1}=p_n+d_n\). By Lemma 4, in Section 3, the inequality\\ \begin{equation}
\left(p_n+d_n\right){}^n=p_{n+1}^n<p_n^{n+(\log  n)^{-1+\epsilon }}
\end{equation}\\
is satisfied for all sufficiently large integer \(n\geq n_0(\epsilon )\). Moreover, taking logarithm of both sides, and simplifying yield\\
\begin{equation}
c_0\frac{n
d_n}{p_n}\leq n \log \left(1+\frac{ d_n}{p_n}\right)<\left((\log  n)^{-1+\epsilon }\right)\log  p_n .
\end{equation}\\
Equivalently, it can be rewritten in the simpler form \\
\begin{equation}
d_n\leq c_1\frac{p_n\log  p_n(\log  n)^{\epsilon }}{n \log  n},
\end{equation}\\
where \(c_0,c_1>0\) are constants. Apply the asymptotic formula for the \textit{ n}th prime \(p_n=(1+o(1))n \log  n\), see (\ref{el31}), and the references
\cite{DP99}, [\cite{NW00}, p. 254], and simplify to obtain \\
\begin{equation}
p_{n+1}-p_n\leq c \log \left(p_n\right){}^{1+\epsilon }
\end{equation}\\
where \(c>0\) is a constant. \\

The result in Theorem 1 is equivalent to the existence of prime numbers in very short intervals \([x,x+y]\) for all sufficiently large real numbers
\(x\geq x_0(\epsilon ) \text{and} y\geq c \log (x)^{1+\epsilon }\), unconditionally. This assertion follows from the inequality\\
\begin{equation}
p_{n+1}-p_n\leq
c \log \left(p_n\right)^{1+\epsilon }\leq c_2 \log (x)^{1+\epsilon },
\end{equation}\\
with $c_2>0$ constant. Here $p_n\leq x$, and $p_{n+1}\leq x+c_2 \log (x)^{1+\epsilon }$. The prime pairs with large prime gaps near
\(\log \left(p_n\right)\log _2\left(p_n\right)\) are studied in \cite{MS07}, and large prime gaps near \(q \log \left(p_n\right)\log _2\left(p_n\right)\)
in arithmetic progressions \(\{ q n+a:\gcd (a,q)=1, n\in \mathbb{N} \}\) are studied in \cite{GP90}. \\

\subsection{Density of Primes in Short Intervals}
The density of primes on large intervals \([x,x+y]\), where $x\geq x_0(\epsilon )$ and $y\geq c x^{1/2+\epsilon }$, has the expected order
of magnitude \\
\begin{equation} \label{el20}
\pi (x+y)-\pi (x)=y/\log  x+o(x/\log  x),
\end{equation}\\
see [\cite{DL12}, p. 70], and \cite{HB88}. In contrast, the precise form of the counting function for primes in short intervals with \(y<c x^{1/2+\epsilon }\) remains
unknown. In fact, the density of prime numbers in very short intervals does not satisfy (\ref{el20}) since it has large fluctuations. More precisely, for
any \(m>1\), the inequalities\\
\begin{equation} 
\log ^{m-1}x<\pi \left(x+\log ^mx\right)-\pi (x)\text{\qquad  and \qquad}\log ^{m-1}x>\pi
\left(x+\log ^mx\right)-\pi (x)
\end{equation} \\

occur infinitely often as \(x\longrightarrow \infty\), tends to infinity, see \cite{MS07}.\\

\textbf{Corollary 2.} Let \(\epsilon >0\) be an arbitrarily small real number, let \(m>1\) be an integer, and let \(x\geq x_0(\epsilon
)\) be a sufficiently large real number. Then, there are at least \\
\begin{equation}
\log ^{m-1-\epsilon }x<\pi \left(x+\log ^{m+\epsilon }x\right)-\pi (x)
\end{equation} \\
prime numbers in the very short interval \(\left[x,x+\log ^{m+\epsilon }x\right]\).\\

A related problem is the resolution of the correct inequality \(\log ^{m-1-\epsilon }x<\pi \left(x+\log ^{m+\epsilon }x\right)-\pi (x)\) \\

for \(x>y>1\) as \(x\longrightarrow \infty\). This is an important problem in the theory of primes numbers, known as the Hardy-Littlewood conjecture.
Equivalently, it asks whether or not the intervals \([y,y+x]\) of lengths \(x\geq 1\) contain more or less primes than the intervals \([1,x]\) as
\(x,y\longrightarrow \infty\) tend to infinity, see [\cite{EP35}, p. 6], [\cite{GR04}, p. 24], \cite{HR74}. \\

The numerical calculations for a potential counterexample is beyond the reach of current machine computations, [\cite{HB88}, p. 217]. For small multiple \(k
x=y\), there is an equality \(\pi (x)=\pi (y+x)-\pi (y)\). \\

\textbf{Maier-Pomerance Integer Gap Problem.} Let \(J(x)=\max _{n\leq x}j(n)\) be the maximal of the Jacobsthal function \(j(n)\). This problem
claims that\\ 
\begin{equation}
J(x)=O\left((\log  x)(\log \log x)^{2+o(1)}\right)
\end{equation} \\
confer [\cite{MP90}, p. 202] and \cite{HT09}. The relevance of this problem arises from the inequality \(J(x)\leq d_{\max }(x)\). Since \(\left.(\log  x)(\text{loglog}
x)^{2+o(1)}\right)\leq c(\log  x)^{1+\epsilon }\) for any constants \(c>0, \text{and} \epsilon >0\), the result in Theorem 1 is in line with this conjecture.\\ 

\textbf{Least Quadratic Nonresidue Problem.} Let \(q\geq 2\), and let \(\chi \neq 1\) be a character modulo \(q\), the least quadratic nonresidue
is defined by \(n_{\chi }=\min \{ n\in \mathbb{N}:\chi (n)\neq 0,1 \}\). Some partial new results for the least quadratic nonresidue for some characters
are demonstrated in [\cite{BG13}, Corollary 2], for example, \(n_{\chi }\ll (\log  q)^{1.37+o(1)}\). This result is strikingly similar to the upper bound stated
in Theorem 1. \\

\section{Supporting Materials}
The basic background information and a few Lemmas are provided in this section. The proofs of these Lemmas use elementary methods. At most a weak form of the Prime Number Theorem \\
\begin{equation} \label{el30}
\pi (x)=\frac{x}{\log x}+o\left (\frac{x}{\log x} \right )
\end{equation} \\
is required, see \cite{RP96}. The innovation here is splicing together these elementary concepts into a coherent and effective result in the theory of primes numbers in very short intervals.\\

\subsection{Formula for the nth Prime, And Associated Series}
The Prime Number Theorem (\ref{el30}) implies that the \textit{n}th prime has an asymptotic expression of the form \(p_n=n \log  n+O( n \text{loglog} n)\).
The Cipolla formula specifies the exact formula for the \textit{n}th prime: \(p_n=(n \log  n)f( n) .\) \\

The function \(f:\mathbb{N}\longrightarrow \mathbb{R}\) has the power series expansion \\ 
\begin{equation} \label{el31}
f(n)=1+\frac{\log\log n-1}{\log  n}+\frac{\log\log
n-2}{\log ^2 n}-\frac{\log ^2\log  n-6\log\log n+11}{\log ^3 n}+\sum _{s\geq 3} \frac{f_s(\log\log n)}{\log ^{s+1} n},
\end{equation}\\
where\\
\begin{equation}
f_s(x)=a_sx^s+a_{s+1}x^{s-1}+\cdots +a_1x+a_0\in \mathbb{Q}[x]
\end{equation}\\
is a polynomial, \(\log ^s n=(\log  n)^s\) is an abbreviated notation for
the \textit{s}th power of the logarithm or its iterated form respectively. For example, \\
\begin{equation}
f_0(x)=x-1,f_1(x)=x-2, 
f_2(x)=-(x^2-6x+11)/2.
\end{equation}\\
The leading coefficient \(a_s=\left.(-1)^{s+1}\right/s\) of the \textit{s}th polynomial \(f_s(x)\) has alternating sign for \(s\geq 1\). Other details
on the asymptotic representation of the \textit{n}th prime \(p_n\) are discussed in [\cite{AT12}, p. 27], \cite{CM02}, and \cite{DP99}.\\

The function \(f(n)\) is quite similar to the simpler function \\
\begin{equation}
g(x)=1+\frac{\text{loglog} x}{\log  x}+\frac{1}{\log  x}\sum _{s\geq 1} (-1)^{s+1}\frac{(\text{loglog}
x)^s}{s \log ^s s}=1+\frac{\text{loglog} x}{\log  x}+\frac{1}{\log  x}\log \left(1+\frac{\text{loglog} x}{\log  x}\right)
\end{equation}\\
and the asymptotic analysis of these two functions are quite similar. \\

An important component in the proof of Lemma 3 is the behavior of the logarithmic ratio\\
\begin{equation}
\log \left ( f(x)/f(x+1) \right ) \text{ as } x \longrightarrow \infty.
\end{equation}\\
The result below works out a lower and an upper estimates of the logarithm ratio of the function \(f(x)\) associated with the \textit{n}th prime \(p_n\).\\

\textbf{Lemma 3.} Let \(n\geq 1\) be a sufficiently large integer. Then \\

(i) $\quad 0<\log  f(n)\leq c_1\frac{\log \log n}
{\log n},$\\

(ii) $\quad 0<f(n)-f(n+1)\leq c_1\frac{\log \log n}{n \log ^2 n},$\\

(iii) $\quad 0<c_2\frac{\log \log n}{n \log
^2 n}\leq \log  \frac{f(n)}{f(n+1)}\leq c_3\frac{\log \log n}{n \log ^2 n},$
\\

where \(c_1,c_2,\text{and} c_3>0\) are constants. \\

\textbf{Proof:} (i) This follows from the power series (\ref{el31}), and the estimate \(\log (1+z)\leq z\), where \(\left| z\right| <1\). \\
(ii) The derivative of the function \(f(x)\) satisfies \(f'(x)<0\) on the interval $[x_0,\infty ) $ for \(x_0\geq 1619\). Hence, the difference \(f(x)-f(x+1)>0\) is nonnegative on the interval
\(\left.\left[x_0,\infty \right.\right)\). \\
(iii) Start with the ratio \\
\begin{equation}
\frac{f(n)}{f(n+1)}=1+\frac{f(n)-f(n+1)}{f(n+1)}.
\end{equation}\\
Taking logarithm and applying the estimate \(\log (1+z)\leq z\), where \(\left| z\right| <1\), lead to \\
\begin{equation}
c_0\frac{f(n)-f(n+1)}{f(n+1)}\leq
\log \left(1+\frac{f(n)-f(n+1)}{f(n+1)}\right)\leq \frac{f(n)-f(n+1)}{f(n+1)}.
\end{equation}\\
where \(c_0>0\) is a constant. In addition, \(1/2\leq f(x)\leq 2\). Consequently, \\
\begin{equation}
c_1(f(n)-f(n+1))\leq \log \left(1+\frac{f(n)-f(n+1)}{f(n+1)}\right)\leq
c_2(f(n)-f(n+1)).
\end{equation}\\
where \(c_1,c_2>0\) are constants. The lower and upper bounds are computed in a term by term basis. For example, the difference of the first term of the power series (\ref{el31}) satisfies \\
\begin{equation}
c_3\frac{\log\log n}{n \log ^2 n}\leq \frac{\log\log n-1}{\log  n}-\frac{\log\log(n+1)-1}{\log
 (n+1)}\leq c_4\frac{\log\log n}{n \log ^2 n},
\end{equation}\\
where \(c_3,c_4>0\) are constants. Lastly, as the function \(f(n)-f(n+1)\) is an alternating power series, the difference of the first term is sufficient.
 $\blacksquare $ \\

The information on the lower and upper bound of the logarithmic ratio \(\log  f(x)/f(x+1)\) provides an effective way of estimating the order of magnitude of the error term in the difference \(p_{n+1}-p_n\). It sidesteps the difficulty encountered while using the asymptotic formula \(p_n=n
\log  n+O( n \log \log n)\). \\

\subsection{The Prime Numbers Equations And Inequalities} 
The exponential prime numbers equations and inequalities \\
\begin{equation}
p_{n+1}^x\pm p_n^y=r(n),\hskip 1 in p_{n+1}^x\pm p_n^y<r(n)
\end{equation}\\
and other variations of these equations and inequalities have been studied by many workers in the field, confer [\cite{GR04}, p. 19], and [\cite{RP96}, p. 256] for discussions and references. The specific case \(p_{n+1}^x\pm p_n^y<0\), which implies the Cramer conjecture, can be established using the elementary method employed in this work. A stronger result considered here seems to be plausible.\\

\textbf{Lemma 4.}  Let \(\epsilon >0\) be an arbitrarily small real number, and let \(n\geq 1\) be a sufficiently large integer. Then \\
\begin{equation}
p_{n+1}^n<p_n^{n+(\log  n)^{-1+\epsilon }}.
\end{equation}\\

\textbf{Proof:} Write a pair of consecutive primes as\\
\begin{equation} \label{el40}
p_n=(n \log  n)f(n)\qquad \text{ and }\qquad p_{n+1}=((n+1)
\log (n+1))f(n+1) .
\end{equation}\\
where \(f(n)\) is the power series (\ref{el31}), and consider the inequality\\
\begin{equation} \label{el41}
\log  p_{n+1}^n<\log  p_n^{n+k},
\end{equation}\\
where \(k=k(n,\epsilon )\) is a function of \(n \text{ and } \epsilon\). Replacing (\ref{el40}) into (\ref{el41}) transforms it into \\
\begin{equation}
\log \left(1+\frac{1}{n}\right)+\log
\left(1+\frac{\log (1+1/n)}{\log  n}\right)+\log  f(n+1)<\frac{k}{n}(\log  n+\text{loglog} n+\log  f(n))+\log  f(n).
\end{equation}\\
Now use \(\log (1+z)\leq z\), where \(\left| z\right| <1\), to reduce it to \\
\begin{equation}
\frac{c_0}{n}+\frac{c_1}{n \log  n}<\frac{k}{n}(\log  n+\log\log
n+\log  f(n))+\log  \frac{f(n)}{f(n+1)},
\end{equation}\\
where \(c_0,c_1>0\) are constants. In Lemma 3 it is shown that \\
\begin{equation}
0<\log  f(n)\leq c_2\frac{\log\log n}{\log  n}\text{                 }\text{  and  }\text{
             }0<\log  \frac{f(n)}{f(n+1)}\leq c_3\frac{\log\log n}{n \log ^2 n},
\end{equation}\\
where \(c_2,c_3>0\) are constants. Accordingly, the penultimate inequality can be written as\\ \begin{equation} \label{el42}
\frac{c_0}{n}+\frac{c_1}{n \log  n}<\frac{k}{n}\left(\log
 n+\log\log n+c_2\frac{\log\log n}{\log  n}\right)+c_3\frac{\log\log n}{n \log ^2 n} .
\end{equation}\\
The first term \(k \log (n)/n\) on the right side of the previous inequality dominates if the parameter is set to \(k=(\log  n)^{-1+\epsilon }\)
for any arbitrarily small real number \(\epsilon >0\). Therefore, the inequality (\ref{el42}) is satisfied by any sufficiently large integer \(n\geq n_0(\epsilon
)\). $\blacksquare $

\section{Applications}
There are many practical and theoretical applications of Theorem 1 in the Mathematics, Cryptography, and Physics. Two obvious applications are considered
here.\\

\subsection{Generating Prime Numbers}
The result sketched in this subsection improves a recent work on the complexity of generating prime numbers in \cite{TH12}.\\

\textbf{Corollary 5.}   Given any real number \(x\in \mathbb{R}\), a prime \(p\geq x\geq 2\) can be constructed in deterministic
polynomial time complexity of \(O\left((\log  x)^8\right)\) arithmetic operations. \\

\textbf{Proof:} Apply the AKS primality test algorithm to the sequence of \(k\geq 1\) consecutive odd numbers\\
\begin{equation}
n+1,n+3,n+5, \ldots, n+k,
\end{equation}\\
where \(k=O\left(\log ^2x\right)\), and \(n=2[ x/2 ]+2\). By Theorem 1, at least one of these integers is a prime number. Thus, repeating the primality
test \(k\) times, result in the running time complexity of \(O\left((\log  x)^8\right)\) arithmetic operations. $\blacksquare $\\

The analysis of the AKS primality test algorithm appears in \cite{AK04}, it has deterministic polynomial time complexity of \(O\left((\log  x)^8\right)\)
arithmetic operations. An improved version of faster running time appears in \cite{LP05}. \\

\subsection{Scherk Expansions of Primes}
The binary expansion of prime number is computable in polynomial time complexity, but the Scherk expansion of a prime has exponential time complexity,
and it seems to be impossible to reduce it to a lower time complexity. This is an expansion of the form, [\cite{NW00}, p. 353], \\
\begin{equation}
p_{2n}=2p_{2n-1}+\sum
_{0\leq k\leq 2n-2} \epsilon _kp_k\text{            }\text{or}\text{         }p_{2n+1}=p_{2n}+\sum _{0\leq k\leq 2n-1} \epsilon _kp_k
\end{equation}\\
where \(\epsilon _k\in \{-1,0,1\}\) and \(p_0=1\). A proof of this expansion is given in the American Mathematics Monthly. For example, the Scherk expansion of the prime
number \\
\begin{equation}
100000000000000000151=(1 \text{or} 2)\cdot 100000000000000000129+\sum _{0\leq k\leq 2n-2} \epsilon _kp_k\text{  }\text{or}\text{  }+\sum
_{0\leq k\leq 2n-1} \epsilon _kp_k
\end{equation}\\
requires about \(2\times 10^{18}\) primes to complete. This amounts to a very large scale computation. However, the short Scherk expansion\\
$$100000000000000000151 = 100000000000000000129+19 + 3$$\\
has polynomial time complexity, this follows from Corollary 5. \\

\section{Numerical Data}
Significant amount of data have been complied by several authors over the past decades, see \cite{NT16}, \cite{NN03}, \cite{PP16}, \cite{TS16}, et alii, for full details. A
few examples of prime gaps, and the estimated upper bound as stated in Theorem 1, using the inequality\\
\begin{equation}
d(n)=15\log \left(p_n\right)\left(\text{loglog}\left(p_n\right)\right)^2-\left(p_{n+1}-p_n\right)
\end{equation}\\
with \(c=15\), see the Maier-Pomerance problem in Section 1 for explanation.\\

\textbf{1.} The largest prime gap for the prime\\
$$p_n= 46242083809774032061673394226721$$\\

listed on the tables in [25], and [26] satisfies the inequality: \\
\begin{equation}
1998=p_{n+1}-p_n\leq 15\log \left(p_n\right)\left(\text{loglog}\left(p_n\right)\right){}^2=4691.026292 .
\end{equation}\\
\textbf{2.} The recently reported very large prime gap of 337446 for the 7976 digits prime \(p_n\), see [40], satisfies the inequality: \\
\begin{equation}
337446=p_{n+1}-p_n\leq 15\log \left(p_n\right)\left(\text{loglog}\left(p_n\right)\right){}^2=2704737.1290 .
\end{equation}\\
\textbf{3.} A few others reported examples for the prime gaps of very large probable primes also satisfy the inequality, in particular,\\
\begin{equation}
2254930=p_{n+1}-p_n\leq 15\log \left(p_n\right)\left(\text{loglog}\left(p_n\right)\right){}^2=36615528.48 .
\end{equation}\\
for a 86853 decimal digits probable prime \(p_n\), see [29], [40].

\end{document}